\documentclass[10pt]{amsart}
\usepackage{amsthm,amsmath,amssymb,tikz}
\usetikzlibrary{arrows}

\title[Near-central permutation factorization]{Near-central Permutation Factorization and Strahov's Generalized Murnaghan-Nakayama Rule}
\author{D. M. Jackson$^1$ and C. A. Sloss$^2$}



\thanks{
${\hspace{-1ex}}^1$ Department of Combinatorics and Optimization, University of Waterloo, Waterloo, Ontario, Canada. Partially supported by an NSERC Discovery Grant. \texttt{dmjackson@math.uwaterloo.ca}}

\thanks{
${\hspace{-1ex}}^2$ Department of Combinatorics and Optimization, University of Waterloo, Waterloo, Ontario, Canada. Partially supported by an NSERC Postgraduate Scholarship. \texttt{csloss@theorem.ca}}

\theoremstyle{plain}
\newtheorem{thm}{Theorem}[section] 
\newtheorem{lemma}[thm]{Lemma} 
 
\newtheorem{cor}[thm]{Corollary} 
\newtheorem{prob}[thm]{Problem}

\newtheorem{defn}[thm]{Definition} 

\theoremstyle{definition}

\renewcommand{\includegraphics}{}

\newcommand{\rj}{\theta}
\newcommand{\nrj}{\vartheta}

\begin{document}

\maketitle

\begin{abstract}
The $(p,q,n)$-dipole problem is a map enumeration problem, arising in perturbative Yang-Mills theory, in which the parameters $p$ and $q$, at each vertex, specify the number of edges separating of two distinguished edges. Combinatorially, it is notable for being a permutation factorization problem which does not lie in the centre of $\mathbb{C}[\mathfrak{S}_n]$, rendering the problem inaccessible through the character theoretic methods often employed to study such problems. This paper gives a solution to this problem on all orientable surfaces when $q=n-1$, which is a combinatorially significant special case: it is a \emph{near-central} problem. We give an encoding of the $(p,n-1,n)$-dipole problem as a product of standard basis elements in the centralizer $Z_1(n)$ of the group algebra $\mathbb{C}[\mathfrak{S}_n]$ with respect to the subgroup $\mathfrak{S}_{n-1}$. 
The generalized characters arising in the solution to the $(p,n-1,n)$-dipole problem are zonal spherical functions of the Gel'fand pair $(\mathfrak{S}_n\times \mathfrak{S}_{n-1}, \mathrm{diag}(\mathfrak{S}_{n-1}))$ and are evaluated explicitly. This solution is used to prove that, for a given surface, the numbers of $(p,n-1,n)$-dipoles and $(n+1-p,n-1,n)$-dipoles are equal, a fact for which we have no combinatorial explanation. These techniques also give a solution to a near-central analogue of the problem of decomposing a full cycle into two factors of specified cycle type. 
\end{abstract}

\section{Introduction}

The character theory of the symmetric group has proven to be a powerful tool for studying combinatorial problems which can be regarded as enumerating the factorizations of a given permutation into factors of specified cycle type. Two major classes of problems which can be studied using such techniques are the enumeration of maps in surfaces (see, for example, \cite{JacksonVisentin:1990}, \cite{JacksonVisentin:2001} and \cite{Jackson:1987}) and the enumeration of ramified covers of the sphere (see, for example, \cite{GouldenJackson:1997} and \cite{GouldenJacksonVakil:2005}). These problems are referred to as \textit{central problems} since they can be encoded as elements of the centre $Z(n)$ of the symmetric group algebra $\mathbb{C}[\mathfrak{S}_n]$. 

Permutation factorization problems involving information other than the cycle type of the factors are an important class of \textit{non-central problems}. Examples of non-central problems include the $(p,q,n)$-dipole problem \cite{ConstableFreedmanHeadrick:2002}, \cite{VisentinWieler:2007}, the non-transitive version of the star factorization problem \cite{JacksonSloss:2011}, and a non-central analogue of the cycle decomposition problem which we introduce in the present paper as Problem \ref{prob:NCCDP}.

The algebraic context we use to study non-central problems is the \textit{centralizer} 
\[
Z_{\mathfrak{H}}(n) := \{ g \in \mathbb{C}[\mathfrak{S}_n]: hgh^{-1}=g \text{ for all } h\in \mathfrak{H}\}
\]
of the group algebra with respect to a subgroup $\mathfrak{H}$. When $\mathfrak{H}= \mathfrak{S}_{n-k}$, we write $Z_{\mathfrak{H}}(n) = Z_k(n)$. Centralizers provide a measure of non-centrality in the sense that, for a combinatorial problem encoded in $Z_{\mathfrak{H}}(n)$, the smaller the subgroup $\mathfrak{H}$, the further the problem departs from centrality. A non-central problem which can be encoded as an element of $Z_1(n)$ is referred to as a \textit{near-central problem}. Combinatorially, a near-central problem is one which depends on the cycle type of a permutation and the length of the cycle containing $n$. In a previous paper \cite{JacksonSloss:2011}, we showed how the solution to a near-central problem can be expressed in terms of \textit{generalized characters} of the symmetric group. (The results we use from this paper are given in Section \ref{sec:Z1-definitions}.) We use Strahov's generalization of the Murnaghan-Nakayama rule \cite{Strahov:2007} to obtain explicit formula for generalized characters corresponding to partitions of the form $(n-k,1)$ and $(n-k-1,2,1)$. These expressions are then applied to a special case of the $(p,q,n)$-dipole problem and the problem of enumerating factorizations of a full cycle. 

\subsection{The loopless dipole problem and a non-central refinement}

\label{sec:loopless-dipoles}

A \textit{dipole} is a 2-cell embedding of a loopless graph, with precisely two vertices, in a locally orientable surface. In the present paper, we consider the case of dipoles in orientable surfaces, where here and throughout we use of the term ``dipole'' to refer to a dipole in an orientable surface. Edge-labelled dipoles may be encoded as a pair of full-cycle permutations $(\sigma_1,\sigma_2)$, corresponding to the edge labels encountered on a counterclockwise circulation of the root and non-root vertices, whose product $\sigma_1\sigma_2$ encodes the face structure of the dipole, in the manner of Kwak and Lee \cite{Kwak:1992}. The genus of the surface can then be determined using the Euler-Poincar\'{e} formula: a dipole with $n$ edges and $m$ faces has genus $\frac12(n-m)$. The problem of enumerating dipoles in a surface of genus $g$ may be solved by computing $K_{(n)}^2$ in $Z(n)$, where $K_{\lambda}$ denotes the sum of all permutations of cycle type $\lambda$, in the following sense: $[\pi]K_{(n)}^2$ is the number of dipoles having face structure corresponding to the permutation $\pi$, which is $[K_{\lambda}]K_{(n)}^2$ whenever $\pi$ has cycle type $\lambda$. In this encoding, $\lambda$ is the partition given by half face-degrees. Here and elsewhere, if $B$ is an element of some basis of an algebra $\mathfrak{A}$, the notation $[B]$ denotes the coefficient extraction operator on $\mathfrak{A}$. Thus, we say that the dipole problem has a \textit{central encoding}. This encoding permits the application of character-theoretic results, namely, those appearing in \cite{Jackson:1987}, to the loopless dipole problem.

It is worth noting that the Kwak-Lee encoding is different than the usual encoding given by the Tutte Embedding Theorem (see \cite{Tutte:1984}, or \cite{JacksonVisentin:2001} for a contemporary account), an encoding in which it is not clear how to forbid loops. From an enumerative point of view, it is important not only to know that a problem is central, but also to have an explicit encoding in terms of standard basis elements of $Z(n)$. 

\subsubsection{The $(p,q,n)$-dipole problem}

A refinement of the dipole problem was introduced by Constable \textit{et al.} \cite{ConstableFreedmanHeadrick:2002} in the study of duality between Yang-Mills theory and string theory. The refined problem deals with dipoles which are \textit{rooted} by the selection of an edge and a vertex. 
\begin{defn}[Root jump and Non-root jump]]
Let $D$ be a rooted dipole with an additional distinguished edge. The \emph{root jump} (resp. \emph{non-root jump}) of $D$, denoted by $\rj(D)$ (resp. $\nrj(D)$), is one plus the number of edges encountered on a counterclockwise circulation of the root (resp. non-root) vertex, starting at the root edge and ending at the second distinguished edge.  (Neither the root edge nor the second distinguished edge are included in the quantities $\rj(D)$ and $\nrj(D)$.) 
\label{defn:root-jump}
\end{defn}

\noindent The weight functions $\rj$ and $\nrj$ are the subject of the following combinatorial question.
\begin{defn}[$(p,q,n)$-dipole problem]
A dipole $D$ with $n$ edges for which $\rj(D)=p$ and $\nrj(D)=q$ is referred to as a \emph{$(p,q,n)$-dipole}. 
The \emph{$(p,q,n)$-dipole problem} is the problem of determining the number $d^{p,q,}_{n,g}$ of $(p,q,n)$-dipoles in an orientable surface of genus $g$. 
\label{defn:pqn-dipole}
\end{defn}

This problem was solved asymptotically on the torus and double torus by Constable \textit{et al.}, and exact solutions for the torus and double torus were given by Visentin and Wieler \cite{VisentinWieler:2007}. The problem is notable from a combinatorial point of view, because for a dipole $D$, the quantities $\rj(D)$ and $\nrj(D)$ cannot be determined from the cycle type of $\pi$ alone, where $\pi$ is the permutation encoding the face structure of $D$. In other words, the problem is non-central.

Kwak and Lee's central encoding for dipoles may be refined in a natural way to give an encoding of a ``labelled analogue'' of the $(p,q,n)$-dipole problem, in which the root edge is labelled $n$, the other distinguished edge is labelled $n-1$, and the other $n-2$ edges are labelled arbitrarily from $\{1,\ldots, n-2\}$. Then the number of labelled $(p,q,n)$-dipoles with face permutation $\pi$ is given by
\begin{equation}
[\pi] \sum_{\substack{\sigma_1\in \mathcal{C}_{(n)} \\ \sigma_1^q(n)=n-1}} \sigma_1 \sum_{\substack{\sigma_2\in \mathcal{C}_{(n)} \\ \sigma_2^p(n)=n-1}} \sigma_2.
\label{eqn:Z2-pqn-encoding}
\end{equation}
Thus, the group algebra element encoding the $(p,q,n)$-dipole problem lies in $Z_2(n)$. Unlike $Z_1(n)$, the algebra $Z_2(n)$ is non-commutative and thus does not have a basis of orthogonal idempotents, which is the main barrier to finding a solution to the $(p,q,n)$-dipole problem through algebraic methods.

In the present paper, we show that when $q=n-1$, there is an encoding for this problem which lies in $Z_1(n)$. (Lemma \ref{lem:p1n-Z1-encoding}.) Thus, the $(p,n-1,n)$-dipole problem is a near-central problem. Near-central enumerative methods are then used to give a solution to the problem on all orientable surfaces. (Theorem \ref{thm:p1n-full-solution}.) As a consequence of this solution, we prove that the $(p,n-1,n)$-dipole problem has an entirely unexpected symmetry in the sense that at present there is no combinatorial proof of it. (Theorem \ref{thm:p1n-symmetry}.)

\subsection{Near-central decompositions of a full cycle}
\label{sec:Z1-decompositions}

Let $C$ be a cycle of length $n$ in $\mathfrak{S}_n$. The problem of determining the number of factorizations of $C$ into factors of specified cycle type has been well studied. The minimal case with two factors was solved by Goulden and Jackson \cite{GouldenJackson:1992} using combinatorial methods. The two-factor case in general was solved by Goupil and Schaeffer \cite{GoupilSchaeffer:1998} and Biane \cite{Biane:2004}. The case with $r$ factors was done by Poulahon and Schaeffer \cite{PoulalhonSchaeffer:2002} and Irving \cite{Irving:2006}, also using character-theoretic means. The central cycle decomposition problem has the following non-central refinement.

\begin{prob}[Near-central Cycle Decomposition Problem]Given $\lambda,\mu\vdash n$ and parts $i$ and $j$ of $\lambda$ and $\mu$ respectively, determine the number of factorizations $\pi_1\pi_2$ of $C$ such that
\begin{enumerate}
\item
$\pi_1$ has cycle type $\lambda$ with $n$ appearing on a cycle of length $i$, and 
\item
$\pi_2$ has cycle type $\mu$ with $n$ appearing on a cycle of length $j$. 
\end{enumerate}
\label{prob:NCCDP}
\end{prob}
An expression for the solution to this problem, as the coefficient of a rational function, is given in Theorem \ref{thm:Z1-cycle-decompositions}. It is natural to explore the near-central version of the problem with the expectation of finding combinatorial structure which is analogous to the combinatorial structure of the central cycle decomposition problem. Understanding such structure would provide an important first step in understanding $Z_2(n)$, which is also a combinatorially interesting algebra. The techniques developed in this paper provide a starting point for such an exploration.

\subsection{Organization of the paper}

This paper is organized as follows. Section \ref{sec:Z1-definitions} presents background definitions and results pertaining to near-central methods in enumerative combinatorics. Section \ref{sec:p1n-encoding} proves that the $(p,n-1,n)$-dipole problem is near-central, and gives an explicit encoding for the problem in $Z_1(n)$. Section \ref{sec:calculating-GC} gives explicit expressions for the generalized characters which arise in the study of $(p,n-1,n)$-dipoles. Finally, Section \ref{sec:combinatorial-applications} uses these expressions to give the generating series for $(p,n-1,n)$-dipoles with respect to number of faces, and an expression for the solution to the near-central cycle decomposition problem. 

\section{Background Definitions and Results}
\label{sec:Z1-definitions}

\subsection{Definitions and Results from Character Theory}

A weakly decreasing sequence of positive integers $\lambda = (\lambda_1,\ldots,\lambda_m)$ is a \emph{partition} of $n$ if $\sum_{1\leq i\leq m} \lambda_i=n$. This is denoted by $\lambda \vdash n$. Each $\lambda_i$ is called a \emph{part} of $\lambda$, and $m_i(\lambda)$ is the number of times $i$ occurs as a part of $\lambda$. Let $m(\lambda)$ denote the number of parts of $\lambda$. The assertion $i\in \lambda$ indicates that $m_i(\lambda)>0$. If $i\in \lambda$, let $i_-(\lambda)$ denote the partition obtained by removing $i$ from $\lambda$ and replacing it with $i-1$. Partitions are often written in the form $(1^{m_1(\lambda)}2^{m_2(\lambda)}\cdots)$. Given $\pi\in \mathfrak{S}_n$, let $\kappa(\pi)$ denote the partition given by the cycle type of $\pi$.

Given a partition $\lambda\vdash n$, the \emph{Ferrers diagram} $\mathcal{F}_{\lambda}$ is a diagram consisting of $n$ boxes arranged in $m(\lambda)$ rows such that, if the parts of $\lambda$ are ordered such that $\lambda_1\geq \lambda_2\geq \cdots \geq \lambda_{m(\lambda)}$, then there are $\lambda_i$ boxes in row $i$, justified to the left margin. A \emph{standard Young tableau} of shape $\lambda$ is a bijective assignment of the integers $\{1,\ldots,n\}$ to the boxes of $\mathcal{F}_{\lambda}$ such that the labels on the boxes increase to the right along rows, and down columns. The set of all standard Young tableaux of shape $\lambda$ is denoted by $\mathrm{SYT}_{\lambda}$. Let $T^*$ be the standard Young tableau obtained by deleting the box with largest label from $T$. Given a standard Young tableau $T$ in which the element $i$ appears in the box in row $j$ and column $k$, let $c_T(i) = k-j$ be the \emph{content} of $i$ in $T$.  It is often convenient to refer to the vector $\mathbf{c}_T = (c_T(1), c_T(2),\ldots, c_T(n))$ as the \emph{content vector} of $T$. In a slight abuse of notation, $\mathbf{c}_{\lambda}$ denotes the multiset of contents of any tableau of shape $\lambda\vdash n$. 

Let $d_{\lambda}$ denote the degree of the irreducible representation indexed by $\lambda$, which is equal to the number of standard Young tableaux of shape $\lambda$. Let $\chi^{\lambda}_{\mu}$ denote the ordinary irreducible character of $\mathfrak{S}_n$ indexed by $\lambda$, evaluated at a permutation of cycle type $\mu$. The following character evaluations are classical, and can be found, for example, in \cite{Jackson:1987}. When $\mu = (n)$, then
\begin{equation}
\chi^{\lambda}_{(n)} = 
\begin{cases}
(-1)^k & \text{~if }\lambda = (n-k,1^k), \\
0 & \text{~otherwise.}
\end{cases}
\label{eqn:char-at-cycle}
\end{equation}
A \emph{hook} is a partitions of the form $(n-k,1^k)$. The evaluation of an irreducible character corresponding to a hook at a permutation of cycle type $\mu$ is given by the formula
\begin{equation}
\chi^{(n-k,1^k)}_{\mu} = [y^k](1+y)^{-1} H_{\mu}(y), \text{ where } H_{\mu}(y) := \prod_{1\leq i\leq m} (1 - (-y)^{\mu_i}).
\label{eqn:hook-gen-series}
\end{equation}

\subsection{The Algebra $Z_1(n)$ and Strahov's generalized characters}

The role played by the centralizer $Z_1(n)$ in the study of the $(p,n-1,n)$-dipole problem is analogous to the role played by the centre of $\mathbb{C}[\mathfrak{S}_n]$ in the study of the loopless dipole problem described in Section \ref{sec:loopless-dipoles}. For $\lambda\vdash n$ and $i\in \lambda$, let $K_{\lambda,i} := \sum_{\pi \in \mathcal{C}_{\lambda,i}} \pi,$ where 
\[
\mathcal{C}_{\lambda,i} := \{\pi \in \mathfrak{S}_n : \kappa(\pi) = \lambda, n \text{ is on a cycle of length } i\}.
\]
The set $\{K_{\lambda,i}\}_{\lambda\vdash n, i\in \lambda}$ forms a linear basis for $Z_1(n)$. 
Let $\mathrm{SYT}_{\lambda,i}$ denote the set of standard Young tableaux of shape $\lambda$ in which the symbol $n$ appears at the end of a row of length $i$. Let $e_T$ denote the Young semi-normal unit corresponding to the tableau $T$. Define the elements
\[
\Gamma^{\lambda,i} := \sum_{T\in \mathrm{SYT}_{\lambda,i}} e_T \in \mathbb{C}[\mathfrak{S}_n].
\]
The following facts about $\Gamma^{\lambda,i}$ are used in this paper and may be found in \cite{JacksonSloss:2011}.

\subsubsection{Generalized characters}
The elements $\Gamma^{\lambda,i}$ lie in $Z_1(n)$ and form a basis of orthogonal idempotents. Thus, for $\lambda,\mu\vdash n$, $i\in \lambda$, $j\in \mu$, the coefficients
\begin{equation}
\gamma^{\lambda,i}_{\mu,j} := \frac{n!}{d_{\lambda}} [K_{\mu,j}] \Gamma^{\lambda,i}
\label{eqn:GC-definition}
\end{equation}
are well-defined. They are called the \emph{generalized characters} of $\mathfrak{S}_n$. 
The definition of generalized characters given by (\ref{eqn:GC-definition}) is equivalent to the definition given by Strahov~\cite{Strahov:2007} as zonal spherical functions of the Gel'fand pair $(\mathfrak{S}_n\times \mathfrak{S}_{n-1}, \mathrm{diag}(\mathfrak{S}_{n-1}))$. 

\subsubsection{Connection coefficients of $Z_1(n)$}
Let $\lambda,\mu,\nu\vdash n$, and let $i\in \lambda$, $j\in \mu$, $k\in \nu$. The \emph{connection coefficients} $c_{\lambda,i,\mu,j}^{\nu,k}$ of $Z_1(n)$ are defined by
\[
K_{\lambda,i} K_{\mu,j} = \sum_{\substack{\nu\vdash n, \\ k\in \nu}} c_{\lambda,i,\mu,j}^{\nu,k} K_{\nu,k}.
\]
They may be expressed in terms of generalized characters as follows:
\begin{equation}
c_{\lambda,i,\mu,j}^{\nu,k} = \frac{|\mathcal{C}_{\lambda,i}||\mathcal{C}_{\mu,j}|}{n!} \sum_{\substack{\rho\vdash n,\\ \ell \in \rho}} \frac{\gamma_{\lambda,i}^{\rho,\ell}\gamma_{\mu,j}^{\rho,\ell} \gamma_{\nu,k}^{\rho,\ell}}{d_{\ell_-(\rho)}} \frac{d_{\rho}}{d_{\ell_-(\rho)}}.
\label{eqn:connection-coefficients}
\end{equation}

\subsubsection{The generalized Diaconis-Greene method}
A generalization of an approach used by Diaconis and Greene \cite{DiaconisGreene:1989} to evaluate ordinary irreducible characters may be refined to give a method to evaluate generalized characters. A polynomial $f\in \mathbb{C}[x_1,\ldots,x_n]$ is said to be \emph{almost symmetric} if, for any $\pi\in \mathfrak{S}_{n-1}$, $f(x_1,\ldots,x_n) = f(x_{\pi(1)}, x_{\pi(2)},\ldots,x_{\pi(n-1)},x_n)$. The ring of almost symmetric polynomials is denoted by $\Lambda^{(1)}[x_1,\ldots,x_n]$. Let 
\[
J_k := \sum_{1\leq i<k} (i,k)
\]
be the $k^{\text{th}}$ \emph{Jucys-Murphy element}. If $f\in \Lambda^{(1)}[x_1,\ldots,x_n]$ is such that $f(J_1,\ldots,J_n) = K_{\lambda,i}$, then for any $\mu\vdash n$ and $j\in \mu$, 
\begin{equation}
\gamma^{\mu,j}_{\lambda,i} = \frac{d_{j_-(\mu)}}{|\mathcal{C}_{\lambda,i}|} f(\mathbf{c}_{j_-(\mu)}, c_{\mu,j}),
\label{eqn:GC-as-content-evaluation}
\end{equation}
where $c_{\mu,j}$ is the content of $n$ in any $T\in \mathrm{SYT}_{\mu,j}$. 

Furthermore, if $\rho\vdash n$, $\ell$ is a part of $\rho$ and $T$ is any standard Young tableau of shape $\rho$, then by the generalized Diaconis-Greene method,
\begin{equation}
\sum_{\substack{\lambda\vdash n, \\ m(\lambda)= m, i\in \lambda}} \frac{|\mathcal{C}_{\lambda,i}|}{d_{\ell_-(\rho)}} \gamma^{\rho,\ell}_{\lambda,i} = [t^{m}] \prod_{1\leq i \leq n} (t+ c_i(T)),
\label{eqn:GC-sum-as-content-polynomials}
\end{equation}
relating a weighted generalized character sum to content polynomials.

\section{A $Z_1(n)$-encoding for the $(p,n-1,n)$-dipole problem}
\label{sec:p1n-encoding}

Although the $(p,q,n)$-dipole problem lies, in general, in $Z_2(n)$, when $q=n-1$ the problem can be encoded in $Z_1(n)$. Working in $Z_1(n)$ instead of $Z_2(n)$ is a significant advantage since the algebra $Z_1(n)$ is well-understood, while the non-commutative algebra $Z_2(n)$ lacks a basis of orthogonal idempotents. The strategy is to encode a dipole not by its pair of vertex permutations as in the encoding given by (\ref{eqn:Z2-pqn-encoding}), but by the vertex permutation $\nu$ corresponding to the non-root vertex and the face permutation $\rho$, isolating those pairs for which $\nu\rho^{-1}$ gives a root vertex permutation corresponding to a specified value of $\theta(D)$. In this approach, it is notationally more convenient to consider dipoles in which the ordinary edges are not labelled, and which can therefore be given a canonical labelling. Namely, specify a canonical cycle $C_p$ for the vertex permutation at the root vertex with the property that $C_p^p(n)=n-1$, say 
\[
C_p = (n, 1,\ldots, p-1, n-1, p,\ldots, n-2).
\]
In the manner of encoding (\ref{eqn:Z2-pqn-encoding}), the number of $(p,n-1,n)$-dipoles with face permutation $\pi$ is given by
\[
D_{\pi,p} := [\pi] \!\!\!\!\!\! \sum_{\substack{\sigma\in \mathcal{C}_{(n)},\\ \sigma^{-1}(n)=n-1}} \!\!\!\!\!\! \sigma C_p.
\]
The observation that
\[
\sum_{\substack{\sigma\in \mathcal{C}_{(n)},\\ \sigma^{-1}(n)=n-1}} \!\!\!\!\!\! \sigma \;= \sum_{\sigma\in \mathcal{C}_{(n-1,1),1}} \!\!\!\!\!\! \sigma (n,n-1) = K_{(n-1,1),1}(n,n-1)
\]
suggests that $D_{\pi,p}$ is ``close'' to being in $Z_1(n)$. The presence of $C_p$ is still an obstruction to carrying out this computation in $Z_1(n)$, but this can be surmounted by noting that
\[
D_{\pi,p} = [C_p^{-1}] \pi^{-1} \sum_{\substack{\sigma\in \mathcal{C}_{(n)},\\ \sigma^{-1}(n)=n-1}} \!\!\!\!\!\! \sigma.
\]
This is the algebraic statement corresponding to the combinatorial observation that encoding dipoles as a (face, vertex) permutation pair is equivalent to encoding them as a (vertex, vertex) permutation pair. Thus, if $S\subset \mathfrak{S}_n$, then the number of $(p,n-1,n)$-dipoles whose face permutation is an element of $S$ is given by
\begin{align*}
D_{\pi,p}  &= \sum_{\pi\in S} [C_p^{-1}]  \pi^{-1}  \sum_{\substack{\sigma\in \mathcal{C}_{(n)},\\ \sigma^{-1}(n)=n-1}} \!\!\!\!\!\! \sigma \\
&= [C_p^{-1}] \left(\sum_{\pi \in S} \pi^{-1}\right) K_{(n-1,1),1} (n,n-1).
\end{align*}
Thus, whenever the set $S$ is invariant under conjugation by $\mathfrak{S}_{n-1}$, the product 
\[
\left(\sum_{\pi \in S} \pi^{-1}\right) K_{(n-1,1),1}
\]
may be computed within $Z_1(n)$, instead of $Z_2(n)$. 

It remains to determine the effect on $Z_1(n)$ basis elements of multiplication by the transposition $(n,n-1)$ and to extract the coefficient of $C_p^{-1}$. If $n$ and $n-1$ are on the same cycle, multiplication by $(n,n-1)$ cuts this cycle in two. Thus, in this case, a full cycle cannot be obtained. If $n$ and $n-1$ are on different cycles, multiplication by $(n,n-1)$ will join them into one cycle. The cycle $C_p^{-1}$ can only be obtained by joining a cycle of length $p$ and a cycle of length $n-p$. Thus,
\begin{align*}
[C_p^{-1}] K_{\lambda,i} (n,n-1) 
&= \begin{cases}
1 & \text{~if } \lambda = (p,n-p) \text{ and } i=p, \\
0 & \text{~otherwise}
\end{cases} \\
&= [K_{(p,n-p),p}] K_{\lambda,i}.
\end{align*}
Combining these facts gives the following encoding for the $(p,n-1,n)$-dipole problem in $Z_1(n)$.
\begin{lemma}
\label{lem:p1n-Z1-encoding}
Let $\lambda\vdash n$, and let $i$ be a part of $\lambda$. Then the number of $(p,n-1,n)$-dipoles having face degree sequence $2\lambda$ in which the root face has degree $2i$ is given by
\[
[K_{(p,n-p),p}] K_{\lambda,i} K_{(n-1,1),1}.
\]
\end{lemma}
We remark that this encoding of the $(p,n-1,n)$-dipole problem provides a means for studying the problem using generalized characters, analogous to the way in which the Kwak-Lee encoding described in Section \ref{sec:loopless-dipoles} allows the application of character theory to the loopless dipole problem. 

\section{Explicit formulae for generalized characters arising in the $(p,n-1,n)$-dipole problem}
\label{sec:calculating-GC}

In light of Lemma \ref{lem:p1n-Z1-encoding} and Equation (\ref{eqn:connection-coefficients}), to solve the $(p,n-1,n)$-dipole problem it suffices to determine an explicit expression for 
\begin{equation}
d^{p,n-1}_{\lambda,i} := \frac{|\mathcal{C}_{\lambda,i}| (n-2)!}{n!} \sum_{\substack{\rho\vdash n,\\ \ell\in \rho}} \frac{\gamma^{\rho,\ell}_{\lambda,i} \gamma^{\rho,\ell}_{(n-1,1),1} \gamma^{\rho,\ell}_{(p,n-p),p}}{d_{\ell_-(\rho)}} \frac{d_{\rho}}{d_{\ell_-(\rho)}}.
\label{eqn:p1n-general-formula}
\end{equation}
This section gives formulae for the generalized characters appearing in this expression.

\subsection{Generalized Characters Evaluated at $((n-1,1),1)$}

\begin{lemma}
Let $\mu\vdash n$ and let $j\in \mu$. Then
\[
\gamma^{\mu,j}_{(n-1,1),1} = \begin{cases}
(-1)^{k} & \text{~if } \mu = (n-k-1,2,1^{k-1}) \text{ and } j=2; \\
(-1)^{k} & \text{~if } \mu = (n-k,1^{k}) \text{ and } j=n-k; \\
(-1)^{k-1} & \text{~if } \mu = (n-k,1^{k}) \text{ and } j=1; \\
0 & \text{~otherwise.}
\end{cases}
\]
\label{lem:GC-evaluation-(n-1)-cycle}
\end{lemma}
\begin{proof}
The standard basis element $K_{(n-1,1),1}$ may be expressed in terms of Jucys-Murphy elements by:
\[
K_{(n-1,1),1} = J_2J_3\cdots J_{n-1}.
\]
Let $T$ be any tableau of shape $\mu$ in which $n$ appears at the end of a row of length $j$. By Equation (\ref{eqn:GC-as-content-evaluation}), 
\[
\gamma^{\mu,j}_{(n-1,1),1} = \frac{d_{j_-(\mu)}}{(n-2)!} \prod_{1\leq k\leq n-1} c_T(k).
\]
This result will be zero when any label other than $n$ occupies a box in the second row and second column. Thus, the only tableaux giving a nonzero result are either hook tableaux or tableaux of shape $(n-k-1,2,1^{k})$ in which $n$ appears at the end of the row of length $2$. The number of tableaux in each case, as well as the product of their contents, are given in Table \ref{table:gamma-n-fixed-calculations}, and they yield the results in the statement of the theorem. 

\begin{table}
\[
\begin{array}{c|cc}
\text{Cycle type } (\mu,j) & d_{j_-(\mu)} &\prod_{2\leq i\leq n-1} c_T(i)  \\
\hline
 (n-k-1,2,1^{k-1}), 2 & \binom{n-2}{k} &(n-k-2)!(-1)^k k!  \\
 (n-k,1^{k}), n-k  &\binom{n-2}{k}   & (n-k-2)! (-1)^k k! \\
 (n-k,1^{k}),1& \binom{n-2}{k-1} & (n-k-1)! (-1)^{k-1}(k-1)! 
\end{array}
\]
\caption{Calculations needed to evaluate $\gamma^{\mu,j}_{(n-1,1),1}$.}
\label{table:gamma-n-fixed-calculations}
\end{table}
\end{proof}

\subsection{Generalized Characters indexed by ``near hook'' partitions}

In light of Lemma \ref{lem:GC-evaluation-(n-1)-cycle}, it suffices to determine the values of the generalized characters $\gamma^{\lambda,i}_{\mu,j}$ when $(\lambda,i)$ is one of $((n-k,1^k),n-k)$, $((n-k,1^k),1)$ or $((n-k-1,2,1^{k-1}),2)$. This may be done using a generalization of the Murnaghan-Nakayama rule due to Strahov \cite{Strahov:2007}, who uses the following terminology. A skew partition $\lambda / \nu$ is called a \emph{broken border strip} if it contains no $2 \times 2$ boxes. (Thus, a broken border strip which is also connected is a rim hook. Two boxes in a Ferrers diagram whose corners touch are not considered to be connected.) A \emph{sharp corner} in a skew diagram is a box which has a box both below it and to the right. A \emph{dull box} has boxes neither to the right nor below it. Let $\mathrm{SC}(\lambda / \nu)$ and $\mathrm{DB}(\lambda / \nu)$ denote the set of sharp corners and dull boxes of $\lambda / \nu$, respectively. Recall that the height of a rim hook $\lambda/\nu$, denoted by $\langle \lambda/\nu \rangle$, is equal to the greatest row occupied by $\lambda/\nu$ minus the least row occupied by $\lambda/\nu$. If $\lambda/\nu$ is a broken border strip, $\langle \lambda/\nu\rangle$ is defined to be the sum of heights of its connected components. Given a skew diagram $\lambda / \nu$ and a part $i$ of $\lambda$, the number $\varphi_{\lambda/\nu,i}$ is defined by
\[
\varphi_{\lambda/\nu,i} = 
(-1)^{\langle \lambda/\nu\rangle} \prod_{s\in \mathrm{SC}(\lambda/\nu)} [c_{\lambda,i} - c(s)] \prod_{\substack{d\in \mathrm{DB}(\lambda/\nu) \\ d\neq \lambda / i_-(\lambda)}}[c_{\lambda,i}- c(d)]^{-1} 
\]
when $\lambda/\nu$ is a broken border strip, and zero otherwise. Strahov's result is as follows.
\begin{thm}[Murnaghan-Nakayama Rule for Generalized Characters --- Strahov \cite{Strahov:2007}]
\label{thm:MN-rule-gen-char}
Let $\lambda,\rho \vdash n$. Let $i$ be a part of $\lambda$, and let $j$ be a part of $\rho$. Then
\[
\gamma^{\lambda,i}_{\mu,j} := \sum_{\substack{\nu\subseteq i_-(\lambda) \\ \nu\vdash n-j }} \varphi_{\lambda/\nu,i} \chi^{\nu}_{\mu \setminus j}.
\]
\end{thm}
The reason this theorem is particularly useful for evaluating $\gamma^{(n-k,1^k),n-k}_{\mu,j}$, $\gamma^{(n-k,1^k),1}_{\mu,j}$ and $\gamma^{(n-k-1,2,1^{k-1}),2}_{\mu,j}$ is that in all three cases, $i_-(\lambda)$ is a hook partition. Thus, every $\nu \subseteq i_-(\lambda)$ is also a hook partition, and $\chi^{\nu}_{\mu\setminus j}$ may be evaluated using Equation (\ref{eqn:hook-gen-series}). The first such result is as follows.
\begin{lemma}
Let $\mu\vdash n$ and let $j$ be a part of $\mu$. Define $R_{n,j}$ by
\[
R_{n,j}(x) := \frac{(n-1) + nx + (-x)^j}{1+x}.
\]
Then, for $0\leq k \leq n-2$,
\begin{equation}
\label{eqn:hook-genchar-n-k}
\gamma^{(n-k,1^k),n-k}_{\mu,j} = \frac{1}{n-1} [x^k] R_{n,j}(x)H_{\mu\setminus j}(x).
\end{equation}
Define $S_{n,j}$ by
\[
S_{n,j}(x): = (-1)^{j-1}\frac{(-1)^jx + nx^j + (n-1)x^{j+1}}{1+x}.
\]
Then, for $1\leq k\leq n-1$,
\begin{equation}
\label{eqn:hook-genchar-1}
\gamma^{(n-k,1^k),1}_{\mu,j} = \frac{1}{n-1}[x^k] S_{n,j}(x) H_{\mu\setminus j}(x).
\end{equation}
\label{lem:hook-genchar-series}
\end{lemma}
\begin{proof}
Equation (\ref{eqn:hook-genchar-n-k}) is proven here. Equation (\ref{eqn:hook-genchar-1}) may be proven in a similar manner. Thus, in the following, $\lambda = (n-k,1^k)$ and $i=n-k$. Furthermore, assume that $k\leq n-k-1$. (The case $k\geq n-k-1$ is similar.)The series $R_{n,j}(x)$ may be expanded as
\[
(n-1) + \sum_{0<\ell<j } (-1)^{\ell+1} x^\ell.
\]
There are three cases:

\noindent \textbf{Case 1}: Suppose $k\geq j$. Then
\begin{align*}
\frac{1}{n-1} [x^k] R_{n,j}(x)H_{\mu\setminus j}(x) &= [x^k]H_{\mu\setminus j}(x) + \sum_{0<\ell<j } \frac{(-1)^{\ell+1}}{n-1} [x^{k-\ell}]H_{\mu\setminus j}(x) \\
&= \chi^{(n-j-k,1^k)}_{\mu\setminus j} + \sum_{0 < \ell < j} \frac{(-1)^{\ell+1}}{n-1} \chi^{(n-j-k+\ell,1^{k-\ell})}_{\mu\setminus j} \\
&= \chi^{(n-j-k,1^k)}_{\mu\setminus j} + \sum_{k-j < \ell < k} \frac{(-1)^{k-\ell+1}}{n-1} \chi^{(n-j-\ell,1^{\ell})}_{\mu\setminus j}.
\end{align*}
In order for $\nu\subseteq i_-(\lambda)$, $\nu$ must be of the form $(n-j-\ell,1^\ell)$ where $k-j+1\leq \ell \leq n-k-1$. The valid range for $\ell$ is illustrated in the following, where the grey boxes are $\nu$, the black box is the distinguished box, and the box labelled ``D'' is a dull box. 
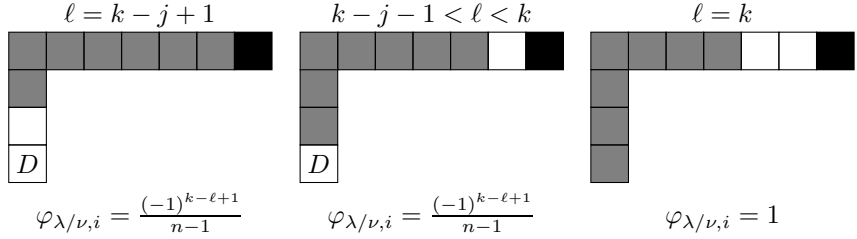
\begin{center}
\begin{tabular}{ccc}
$\ell = k-j+1$ & $k-j-1<\ell < k$ & $\ell = k$ \\

\begin{tikzpicture}[scale=0.5]


\fill[color=gray] (0,0) rectangle (6,-1); 
\fill[color=black] (6,0) rectangle (7,-1);

\draw (0,0) rectangle (1,-1);
\draw (1,0) rectangle (2,-1);
\draw (2,0) rectangle (3,-1);
\draw (3,0) rectangle (4,-1);
\draw (4,0) rectangle (5,-1);
\draw (5,0) rectangle (6,-1);
\draw (6,0) rectangle (7,-1);


\fill[color=gray] (0,-1) rectangle (1,-2); 

\draw (0,-1) rectangle (1,-2);
\draw (0,-2) rectangle (1,-3);
\draw (0,-3) rectangle (1,-4);

\path (0.5,-3.5) node[draw = none, fill=none] (n1) {$D$};

\end{tikzpicture}

&

\begin{tikzpicture}[scale=0.5]


\fill[color=gray] (0,0) rectangle (5,-1); 
\fill[color=black] (6,0) rectangle (7,-1);

\draw (0,0) rectangle (1,-1);
\draw (1,0) rectangle (2,-1);
\draw (2,0) rectangle (3,-1);
\draw (3,0) rectangle (4,-1);
\draw (4,0) rectangle (5,-1);
\draw (5,0) rectangle (6,-1);
\draw (6,0) rectangle (7,-1);


\fill[color=gray] (0,-1) rectangle (1,-3); 

\draw (0,-1) rectangle (1,-2);
\draw (0,-2) rectangle (1,-3);
\draw (0,-3) rectangle (1,-4);

\path (0.5,-3.5) node[draw = none, fill=none] (n1) {$D$};

\end{tikzpicture}

&

\begin{tikzpicture}[scale=0.5]


\fill[color=gray] (0,0) rectangle (4,-1); 
\fill[color=black] (6,0) rectangle (7,-1);

\draw (0,0) rectangle (1,-1);
\draw (1,0) rectangle (2,-1);
\draw (2,0) rectangle (3,-1);
\draw (3,0) rectangle (4,-1);
\draw (4,0) rectangle (5,-1);
\draw (5,0) rectangle (6,-1);
\draw (6,0) rectangle (7,-1);


\fill[color=gray] (0,-1) rectangle (1,-4); 

\draw (0,-1) rectangle (1,-2);
\draw (0,-2) rectangle (1,-3);
\draw (0,-3) rectangle (1,-4);

\end{tikzpicture}

\\
$\varphi_{\lambda/\nu,i} = \frac{(-1)^{k-\ell+1}}{n-1}$ & $\varphi_{\lambda/\nu,i} =\frac{(-1)^{k-\ell+1}}{n-1}$ & $\varphi_{\lambda/\nu,i} = 1$
\end{tabular}
\end{center}
\vspace{0.5cm}
\noindent
Applying Theorem \ref{thm:MN-rule-gen-char},
\begin{align*}
\gamma^{(n-k,1^k),n-k}_{\mu,j} &= \chi^{(n-j-k,1^k)}_{\mu\setminus j} + \sum_{k-j < \ell < k} \frac{(-1)^{k-\ell+1}}{n-1} \chi^{(n-j-\ell,1^{\ell})}_{\mu\setminus j} \\
&= \frac{1}{n-1} [x^k] R_{n,j}(x)H_{\mu\setminus j}(x).
\end{align*}

\noindent \textbf{Case 2}: Suppose $k<(n-j)<n-k$.
In this case, as well as in Case 3, the range $0<\ell<j$ of the index of summation is truncated since the summand is zero for certain values of $\ell$. Since $j>k$, 
\begin{align*}
\frac{1}{n-1} [x^k] R_{n,j}(x)H_{\mu\setminus j}(x) &= [x^k]H_{\mu\setminus j}(x) + \sum_{0<\ell<j } \frac{(-1)^{\ell+1}}{n-1} [x^{k}]x^{\ell} H_{\mu\setminus j}(x)  \\
&= [x^k]H_{\mu\setminus j}(x) + \sum_{0<\ell\leq k} \frac{(-1)^{\ell+1}}{n-1} [x^{k-\ell}]H_{\mu\setminus j}(x) \\
&= \chi^{(n-j-k,1^k)}_{\mu\setminus j} + \sum_{0< \ell \leq k} \frac{(-1)^{\ell+1}}{n-1} \chi^{(n-j-k+\ell,1^{k-\ell})}_{\mu\setminus j} \\
&= \chi^{(n-j-k,1^k)}_{\mu\setminus j} + \sum_{0\leq \ell < k} \frac{(-1)^{k- \ell+1}}{n-1} \chi^{(n-j-\ell,1^{\ell})}_{\mu\setminus j}.
\end{align*}
In this case, $\nu = (n-j-\ell,1^{\ell})$ where $0\leq \ell \leq k$. The range of validity for $\ell$ is illustrated as follows.
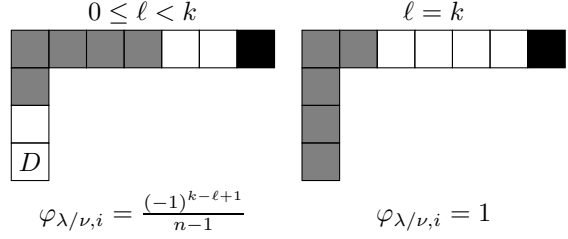
\begin{center}
\begin{tabular}{cc}
$0\leq \ell < k$ & $\ell=k$ \\

\begin{tikzpicture}[scale=0.5]


\fill[color=gray] (0,0) rectangle (4,-1); 
\fill[color=black] (6,0) rectangle (7,-1);
\fill[color=gray] (0,-1) rectangle (1,-2);

\draw (0,0) rectangle (1,-1);
\draw (1,0) rectangle (2,-1);
\draw (2,0) rectangle (3,-1);
\draw (3,0) rectangle (4,-1);
\draw (4,0) rectangle (5,-1);
\draw (5,0) rectangle (6,-1);
\draw (6,0) rectangle (7,-1);



\draw (0,-1) rectangle (1,-2);
\draw (0,-2) rectangle (1,-3);
\draw (0,-3) rectangle (1,-4);

\path (0.5,-3.5) node[draw = none, fill=none] (n1) {$D$};

\end{tikzpicture}

&

\begin{tikzpicture}[scale=0.5]


\fill[color=gray] (0,0) rectangle (2,-1); 
\fill[color=black] (6,0) rectangle (7,-1);
\fill[color=gray] (0,-1) rectangle (1,-4); 

\draw (0,0) rectangle (1,-1);
\draw (1,0) rectangle (2,-1);
\draw (2,0) rectangle (3,-1);
\draw (3,0) rectangle (4,-1);
\draw (4,0) rectangle (5,-1);
\draw (5,0) rectangle (6,-1);
\draw (6,0) rectangle (7,-1);


\draw (0,-1) rectangle (1,-2);
\draw (0,-2) rectangle (1,-3);
\draw (0,-3) rectangle (1,-4);

\end{tikzpicture}

\\
$\varphi_{\lambda/\nu,i} = \frac{(-1)^{k-\ell+1}}{n-1}$ & $\varphi_{\lambda/\nu,i} = 1$
\end{tabular}
\end{center}

\vspace{0.5cm} \noindent
By Theorem \ref{thm:MN-rule-gen-char},
\begin{align*}
\gamma^{(n-k,1^k),n-k}_{\mu,j} &= \chi^{(n-j-k,1^k)}_{\mu\setminus j} + \sum_{0\leq \ell < k} \frac{(-1)^{k- \ell+1}}{n-1} \chi^{(n-j-\ell,1^{\ell})}_{\mu\setminus j} \\
&= \frac{1}{n-1} [x^k] R_{n,j}(x)H_{\mu\setminus j}(x).
\end{align*}

\noindent \textbf{Case 3}: Suppose $n\leq k+j$.
This case relies on the fact that $H_{\mu\setminus j}(x)$ is a polynomial of degree $n-j-1$, so $[x^i]H_{\mu\setminus j}(x)=0$ when $i\geq n-j$. In particular, $[x^k]H_{\mu\setminus j}(x)=0$. Furthermore, $j\geq n-k>k$. Thus, in this case,
\begin{align*}
\frac{1}{n-1} [x^k] R_{n,j}(x)H_{\mu\setminus j}(x) &= [x^k]H_{\mu\setminus j}(x) + \sum_{0<\ell <j} \frac{(-1)^{\ell+1}}{n-1} [x^{k}]x^{\ell} H_{\mu\setminus j}(x) \\
&= \sum_{k-n+j+1\leq \ell \leq k} \frac{(-1)^{\ell+1}}{n-1} [x^{k-\ell}] H_{\mu\setminus j}(x) \\
&= \sum_{k-n+j+1\leq \ell \leq k} \frac{(-1)^{\ell+1}}{n-1} \chi^{(n-j-k+\ell,1^{k-\ell})}_{\mu,j} \\
&= \sum_{0 \leq \ell \leq n-j-1} \frac{(-1)^{k-\ell+1}}{n-1} \chi^{(n-j-\ell,1^{\ell})}_{\mu,j}.
\end{align*}
Since $n-j\leq k\leq n-k-1$, then all $\nu = (n-j-\ell, 1^{\ell})$ such that $0\leq \ell\leq n-j-1$ satisfy $\nu\subseteq i_-(\lambda)$. In all such cases, $\varphi_{\lambda/\nu,i} = \frac{(-1)^{k-\ell+1}}{n-1}$. 
Thus, by Theorem \ref{thm:MN-rule-gen-char},
\begin{align*}
\gamma^{(n-k,1),n-k}_{\mu,j} &= \sum_{0 \leq \ell \leq n-j-1} \frac{(-1)^{k-\ell+1}}{n-1} \chi^{(n-j-\ell,1^{\ell})}_{\mu,j} \\
&= \frac{1}{n-1} [x^k] R_{n,j}(x)H_{\mu\setminus j}(x).
\end{align*}

Thus, when $k\leq n-k-1$, for all values of $j$,
\[
\frac{1}{n-1} [x^k] R_{n,j}(x)H_{\mu\setminus j}(x) = \gamma^{(n-k,1^k),n-k}_{\mu,j}.
\]
If $k\geq n-k-1$, the result may be proven in a similar manner.
\end{proof}

\subsubsection{Specializations}
Specializing this to the generalized characters arising in the $(p,n-1,n)$-dipole problem gives the following explicit evaluations.
\begin{cor}
\label{cor:gen-char-two-part-evaluations-A}
Let $0\leq k\leq n-2$ and $1\leq p\leq n-1$. 
If $k\leq n-k-1$, then
\[
\gamma^{(n-k,1^k),n-k}_{(n-p,p),p} = 
\begin{cases}
(-1)^{k-1} \frac{n-p}{n-1} & \text{~if } (n-p)\leq k, \\
(-1)^k \frac{n-k-1}{n-1} &\text{~if } k<(n-p) <n-k, \\
(-1)^k \frac{n-p}{n-1} & \text{~if } n-k\leq (n-p).
\end{cases}
\]
If $k\geq n-k-1$, then
\[
\gamma^{(n-k,1^k),n-k}_{(n-p,p),p} = 
\begin{cases}
(-1)^{k-1} \frac{n-p}{n-1} &\text{~if } (n-p) < n-k, \\
(-1)^{k-1} \frac{n-k-1}{n-1} & \text{~if } n-k\leq (n-p) \leq k, \\
(-1)^k \frac{n-p}{n-1} & \text{~if } k< (n-p).
\end{cases}
\]
\end{cor}

\begin{cor}
Let $1\leq k\leq n-1$, and let $1\leq p \leq n-1$. If $k\leq n-k-1$, then 
\[
\gamma^{(n-k,1^k),1}_{(n-p,p),p} = 
\begin{cases}
(-1)^{k} \frac{n-p}{n-1} & \text{~if } (n-p)\leq k, \\
(-1)^k \frac{k}{n-1} &\text{~if } k<(n-p) <n-k, \\
(-1)^{k-1} \frac{n-p}{n-1} & \text{~if } n-k\leq (n-p).
\end{cases}
\]
If $k\geq n-k-1$, then 
\[
\gamma^{(n-k,1^k),1}_{(n-p,p),p} = 
\begin{cases}
(-1)^{k} \frac{n-p}{n-1} &\text{~if } (n-p) < n-k, \\
(-1)^{k-1} \frac{k}{n-1} & \text{~if } n-k\leq (n-p) \leq k, \\
(-1)^{k-1} \frac{n-p}{n-1} & \text{~if } k< (n-p).
\end{cases}
\]
\label{lem:gen-char-two-part-evaluations-B}
\end{cor}

In the case of $\gamma^{(n-k-1,2,1^{k-1}),2}_{(n-p,p),p}$, although polynomials analogous to $R_{n,j}$ and $S_{n,j}$ are not known, it is nevertheless possible to evaluate generalized characters of this form, as follows.
\begin{lemma}
Let $1\leq k\leq n-3$, and let $2\leq p\leq n-1$. If $k\leq n-k-2$, then
\[
\gamma^{(n-k-1,2,1^{k-1}),2}_{(n-p,p),p} = 
\begin{cases}
(-1)^k \frac{n-p}{k(n-k-2)} & \text{~if } (n-p)\leq k, \\
0 & \text{~if } k<(n-p)<n-k-1, \\
(-1)^{k+1}\frac{n-p}{k(n-k-2)} & \text{~if } (n-p)\geq n-k-1.
\end{cases}
\]
If $k\geq n-k-2$, then
\[
\gamma^{(n-k-1,2,1^{k-1}),2}_{(n-p,p),p} = 
\begin{cases}
(-1)^k \frac{n-p}{k(n-k-2)} & \text{~if } (n-p)\leq n-k-2, \\
0 & \text{~if } n-k-2<(n-p)\leq k, \\
(-1)^{k+1}\frac{n-p}{k(n-k-2)} & \text{~if } (n-p)>k.
\end{cases}
\]
\label{lem:gen-char-two-part-evaluations-C}
\end{lemma}
\begin{proof}
The proof is by a case analysis similar to the proofs of Lemma \ref{lem:hook-genchar-series}. Throughout the following, $\lambda = (n-k-1,2,1^{k-1})$. Each case is illustrated with a diagram in which grey boxes indicate $\nu$, the black box indicates the distinguished box, and white boxes indicate $2_-(\lambda) /\nu$. Sharp corners and dull boxes of $\lambda / \nu$ are indicated on diagrams by $S$ and $D$, respectively, with the exception of the distinguished box, which is always a dull box. Throughout the following, the fact that
\[
\chi_{(n-p)}^{n-p-\ell,1^{\ell}} = (-1)^{\ell}
\]
is used. 

\noindent
\textbf{Boundary cases}: the cases $k=1$ and $k=n-3$ are treated differently than $2\leq k\leq n-4$, but give the same result. Suppose $k=1$. (The $k=n-3$ case is very similar.) If $p=n-1$, then $(n-p)\leq k$, and a typical diagram has the following form. 

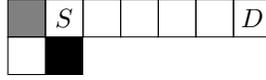
\begin{center}
\begin{tikzpicture}[scale=0.5]


\fill[color=gray] (0,0) rectangle (1,-1); 

\draw (0,0) rectangle (1,-1);
\draw (1,0) rectangle (2,-1);
\draw (2,0) rectangle (3,-1);
\draw (3,0) rectangle (4,-1);
\draw (4,0) rectangle (5,-1);
\draw (5,0) rectangle (6,-1);
\draw (6,0) rectangle (7,-1);


\draw (0,-1) rectangle (1,-2);

\path (1.5,-0.5) node[draw = none, fill=none] (n1) {$S$};
\path (6.5,-0.5) node[draw = none, fill=none] (n1) {$D$};


\fill[color=black] (1,-1) rectangle (2,-2);

\end{tikzpicture}
\end{center}

In this case, the formula in the statement of the Lemma gives
\[
\gamma^{(n-2,2),2}_{(n-1,1),n-1} = \frac{(-1)^k}{k(n-k-2)} = -\frac{1}{n-3}.
\]
This agrees with the value given by Theorem \ref{thm:MN-rule-gen-char}, namely
\[
\gamma^{(n-2,2),2}_{(n-1,1),n-1} = -\frac{1}{n-3}\chi_{(1)}^{(1)} = -\frac{1}{n-3}.
\]
When $3\leq p\leq n-2$, then $k<(n-p)<(n-k-1)$. In this case, the only partitions $\nu\vdash n-j$ which are contained in $(n-2,1)$ are $(n-p)$ and $(n-p-1,1)$:

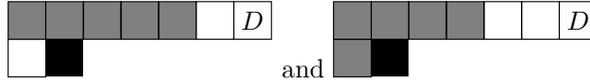
\begin{center}
\begin{tikzpicture}[scale=0.5]


\fill[color=gray] (0,0) rectangle (5,-1); 

\draw (0,0) rectangle (1,-1);
\draw (1,0) rectangle (2,-1);
\draw (2,0) rectangle (3,-1);
\draw (3,0) rectangle (4,-1);
\draw (4,0) rectangle (5,-1);
\draw (5,0) rectangle (6,-1);
\draw (6,0) rectangle (7,-1);


\draw (0,-1) rectangle (1,-2);


\fill[color=black] (1,-1) rectangle (2,-2);

\path (6.5,-0.5) node[draw = none, fill=none] (n1) {$D$};

\end{tikzpicture}
and
\begin{tikzpicture}[scale=0.5]


\fill[color=gray] (0,0) rectangle (4,-1); 
\fill[color=gray] (0,-1) rectangle (1,-2); 

\draw (0,0) rectangle (1,-1);
\draw (1,0) rectangle (2,-1);
\draw (2,0) rectangle (3,-1);
\draw (3,0) rectangle (4,-1);
\draw (4,0) rectangle (5,-1);
\draw (5,0) rectangle (6,-1);
\draw (6,0) rectangle (7,-1);


\draw (0,-1) rectangle (1,-2);


\fill[color=black] (1,-1) rectangle (2,-2);

\path (6.5,-0.5) node[draw = none, fill=none] (n1) {$D$};

\end{tikzpicture}
.
\end{center}

In both cases, $\lambda/\nu$ has height 0, no sharp corners, and one dull box (of content $n-3$) aside from the dull box corresponding to the distinguished part 2, so
\[
\varphi_{\lambda /\nu, 2} = \frac{1}{0 - (n-k-2)} =  -\frac{1}{n-3}.
\]
Thus, Theorem \ref{thm:MN-rule-gen-char} gives
\[
\gamma^{(n-2,2),2}_{(n-p,p),p} = -\frac{1}{n-3}(\chi^{(n-p)}_{(n-p)} + \chi^{(n-p-1,1)}_{(n-p)}) = 0.
\]
When $p=2$, $(n-p)\geq (n-k-1)$, and the diagrams corresponding to the two possibilities for $\nu$ are 
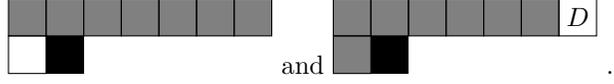
\begin{center}
\begin{tikzpicture}[scale=0.5]


\fill[color=gray] (0,0) rectangle (7,-1); 

\draw (0,0) rectangle (1,-1);
\draw (1,0) rectangle (2,-1);
\draw (2,0) rectangle (3,-1);
\draw (3,0) rectangle (4,-1);
\draw (4,0) rectangle (5,-1);
\draw (5,0) rectangle (6,-1);
\draw (6,0) rectangle (7,-1);


\draw (0,-1) rectangle (1,-2);


\fill[color=black] (1,-1) rectangle (2,-2);

\end{tikzpicture}
and
\begin{tikzpicture}[scale=0.5]


\fill[color=gray] (0,0) rectangle (6,-1); 
\fill[color=gray] (0,-1) rectangle (1,-2); 

\draw (0,0) rectangle (1,-1);
\draw (1,0) rectangle (2,-1);
\draw (2,0) rectangle (3,-1);
\draw (3,0) rectangle (4,-1);
\draw (4,0) rectangle (5,-1);
\draw (5,0) rectangle (6,-1);
\draw (6,0) rectangle (7,-1);


\draw (0,-1) rectangle (1,-2);


\fill[color=black] (1,-1) rectangle (2,-2);

\path (6.5,-0.5) node[draw = none, fill=none] (n1) {$D$};

\end{tikzpicture}
.
\end{center}
The argument proceeds as in the case when $3\leq p\leq n-2$, with the exception that when $\nu = (n-p)$, $\lambda / \nu$ has no dull boxes, so $\varphi_{\lambda/ \nu, 2} = 1$. Thus,
\[
\gamma^{(n-2,2),2}_{(n-2,2),2} = \chi^{(n-2)}_{(n-2)} - \frac{1}{n-3} \chi^{(n-3,1)}_{(n-2)} = \frac{n-2}{n-3}.
\]

\noindent 
\textbf{General case}: for $k\leq 2 \leq n-4$, there are six subcases to consider. As they are similar to each other, one is singled out here for a detailed presentation. This case has been selected because it illustrates all the peculiarities which must be taken into account when computing $\gamma^{(n-k-1,2,1^{k-1}),2}_{(n-p,p),p}$, but did not arise when computing $\gamma^{(n-k,1^k),n-k}_{\mu,j}$ and $\gamma^{(n-k,1^k),1}_{\mu,j}$. Suppose $k\geq n-k-2$ and $(n-k-2)<(n-p)\leq k$. Let $\nu = (n-p-\ell,1^{\ell})$. The range of $\ell$ for which $\nu \subseteq 2_-(\lambda)$ is $k-p+1\leq  \ell \leq n-p-1$. The values of $\varphi_{\lambda/\nu,2}$ for values of $\ell$ in this range are determined based on the following three cases:

\begin{center}
\begin{tabular}{ccc}
\begin{tikzpicture}[scale=0.5]


\fill[color=gray] (0,0) rectangle (4,-1); 
\fill[color=gray] (0,-1) rectangle (1,-3); 

\draw (0,0) rectangle (1,-1);
\draw (1,0) rectangle (2,-1);
\draw (2,0) rectangle (3,-1);
\draw (3,0) rectangle (4,-1);


\draw (0,-1) rectangle (1,-2);
\draw (0,-2) rectangle (1,-3);
\draw (0,-3) rectangle (1,-4);
\draw (0,-4) rectangle (1,-5);
\draw (0,-5) rectangle (1,-6);
\draw (0,-6) rectangle (1,-7);


\fill[color=black] (1,-1) rectangle (2,-2);

\path (0.5,-6.5) node[draw = none, fill=none] (n1) {$D$};

\end{tikzpicture}
&
\begin{tikzpicture}[scale=0.5]


\fill[color=gray] (0,0) rectangle (3,-1); 
\fill[color=gray] (0,0) rectangle (1,-4); 

\draw (0,0) rectangle (1,-1);
\draw (1,0) rectangle (2,-1);
\draw (2,0) rectangle (3,-1);
\draw (3,0) rectangle (4,-1);


\draw (0,-1) rectangle (1,-2);
\draw (0,-2) rectangle (1,-3);
\draw (0,-3) rectangle (1,-4);
\draw (0,-4) rectangle (1,-5);
\draw (0,-5) rectangle (1,-6);
\draw (0,-6) rectangle (1,-7);


\fill[color=black] (1,-1) rectangle (2,-2);

\path (0.5,-6.5) node[draw = none, fill=none] (n1) {$D$};
\path (3.5,-0.5) node[draw = none, fill=none] (n2) {$D$};

\end{tikzpicture}
&
\begin{tikzpicture}[scale=0.5]


\fill[color=gray] (0,0) rectangle (1,-6); 

\draw (0,0) rectangle (1,-1);
\draw (1,0) rectangle (2,-1);
\draw (2,0) rectangle (3,-1);
\draw (3,0) rectangle (4,-1);


\draw (0,-1) rectangle (1,-2);
\draw (0,-2) rectangle (1,-3);
\draw (0,-3) rectangle (1,-4);
\draw (0,-4) rectangle (1,-5);
\draw (0,-5) rectangle (1,-6);
\draw (0,-6) rectangle (1,-7);


\fill[color=black] (1,-1) rectangle (2,-2);

\path (0.5,-6.5) node[draw = none, fill=none] (n1) {$D$};
\path (3.5,-0.5) node[draw = none, fill=none] (n2) {$D$};
\path (1.5,-0.5) node[draw = none, fill=none] (n3) {$S$};

\end{tikzpicture}

\\
(1) $\ell = k-p+1$ &\quad (2) $k-p+1<\ell < n-p-1$\qquad & (3) $\ell = n-p-1$ 
\end{tabular}
\end{center}

\begin{enumerate}
\item
When $\ell = k-p+1$, $\nu = (n-k-1,1^{k-j+1})$, so $\lambda / \nu$ has height $j$, no sharp corners, and one dull box (of content $-k$) aside from the distinguished dull box. Thus, in this case,
\[
\varphi_{\lambda/\nu, 2} = \frac{(-1)^p}{k}.
\]

\item
When $k-p+1<\ell < n-p-1$, $\lambda/\nu$ has height $k-\ell-1$, no sharp corners, and two dull boxes aside from the distinguished box, having contents $(n-k-2)$ and $-k$. Thus,
\[
\varphi_{\lambda/\nu, 2} = \frac{(-1)^{k-\ell-1}}{(0 - (-k))(0 - (n-k-2)} =  \frac{(-1)^{k-\ell}}{k(n-k-2)}.
\]

\item
Finally, when $\ell = n-p-1$, $\nu = (1^{n-p})$, so $\lambda / \nu$ has height $k-\ell$, one sharp corner of content $1$, and two dull boxes aside from the distinguished box, having contents $(n-k-2)$ and $-k$. Thus,
\[
\varphi_{\lambda/\nu, 2} = \frac{(-1)^{k-\ell} (0-1)}{0-(-k))(0 - (n-k-2)} = \frac{(-1)^{k-\ell}}{k(n-k-2)}.
\]
\end{enumerate}
Thus, by Theorem \ref{thm:MN-rule-gen-char},
\begin{align*}
\gamma^{(n-k-1,2,1^{k-1}),2}_{(n-p,p),p} &= \frac{(-1)^p}{k} \chi_{(n-p)}^{(n-k-1,1^{k-p+1})} + \sum_{k-p-1<\ell<n-p-1} \frac{(-1)^{k-\ell}}{k(n-k-2)} \chi_{(n-p)}^{(n-p-\ell,1^{\ell})} \\ 
&\qquad + \frac{(-1)^{k-n+p+1}}{k(n-k-2)} \chi_{(n-p)}{(1^{n-p})} \\
&= \frac{(-1)^{k+1}}{k} + \frac{(-1)^k (n-k-1)}{k(n-k-2)} + \frac{(-1)^{k+1}}{k(n-k-2)}=0,
\end{align*}
completing the proof.
\end{proof}

Having determined explicit formulae for the generalized characters arising in the study of the $(p,n-1,n)$-dipole problem, we may now turn our attention to combinatorial applications of these results.

\section{Combinatorial Applications}
\label{sec:combinatorial-applications}

\subsection{Generating series for $(p,n-1,n)$-dipoles}

The generalized character formulae given in Section \ref{sec:calculating-GC} may be used to express the solution to the $(p,n-1,n)$-dipole problem in terms of generalized characters, as follows. 
\begin{thm}
\label{thm:p1n-dipoles-general-formula}
Let $\lambda\vdash n$ and let $i\in \lambda$. Let $1\leq p \leq n-1$. Then the number of $(p,n-1,n)$-dipoles (with unlabelled ordinary edges) having face degree sequence $2\lambda$ and a root face of degree $2i$ is given by
\[
d^{p,n-1}_{\lambda,i} = \frac{|\mathcal{C}_{\lambda,i}| (n-2)!}{n!} \left(A^{\lambda,i}_{n,p} + B^{\lambda,i}_{n,p} + C^{\lambda,i}_{n,p}\right),
\]
where
\[
A^{\lambda,i}_{n,p} = \sum_{0\leq k\leq n-2} \frac{(-1)^k\gamma^{(n-k,1^k),n-k}_{\lambda,i}\gamma^{(n-k,1^k),n-k}_{(p,n-p),p}}{d_{(n-k-1,1^k)}} \frac{n-1}{n-k-1},
\]
\[
B^{\lambda,i}_{n,p} = \sum_{1\leq k \leq n-1} \frac{(-1)^{k-1} \gamma^{(n-k,1^k),1}_{\lambda,i} \gamma^{(n-k,1^k),1}_{(p,n-p),p}}{d_{(n-k,1^{k-1})}} \frac{n-1}{k},
\]
and
\[
C^{\lambda,i}_{n,p} = \sum_{1\leq k \leq n-3} \frac{(-1)^k \gamma^{(n-k-1,2,1^{k-1}),2}_{\lambda,i} \gamma^{(n-k-1,2,1^{k-1}),2}_{(p,n-p),p}}{d_{(n-k-1,1^k)}}\frac{nk(n-k-2)}{(n-k-1)(k+1)}.
\]
\end{thm}
\begin{proof}
This result follows from Equation (\ref{eqn:p1n-general-formula}) and Lemma \ref{lem:GC-evaluation-(n-1)-cycle}.
\end{proof}
This result becomes more explicit when summing over all partitions corresponding to a given surface by using the two following lemmas.
\begin{lemma}
Let $\rho = (n-k,1^k)$, and let $\ell \in \{1,n-k\}$. Then
\[
\sum_{\substack{\lambda\vdash n, \\ m(\lambda)= m, i\in \lambda}} \frac{|\mathcal{C}_{\lambda,i}|}{d_{\ell_-(\rho)}} \gamma^{\rho,\ell}_{\lambda,i} = [t^m] n! \binom{t+n-k-1}{n}.
\]
\label{lem:GC-sum-A}
\end{lemma}
\begin{proof}
This follows from Equation (\ref{eqn:GC-sum-as-content-polynomials}). The contents along the first row of a tableau of shape $(n-k,1^k)$ are $0,1,2,\ldots,n~-~k~-~1$, and the contents in the first column (excluding the box in the first row) are $-1,-2,\ldots,-k$. Thus,
\[
\prod_{1\leq i \leq n} (t+ c_i(T)) = \prod_{-k \leq i \leq n-k-1} (t+i),
\]
from which the result follows.
\end{proof}
\begin{lemma}
Let $\rho = (n-k-1,2,1^{k-1})$ and let $\ell\in \{n-k-1,2,1\}$. Then
\[
\sum_{\substack{\lambda\vdash n, \\ m(\lambda)= m, i\in \lambda}} \frac{|\mathcal{C}_{\lambda,i}|}{d_{\ell_-(\rho)}} \gamma^{\rho,\ell}_{\lambda,i} = [t^m] (n-1)! t \binom{t+n-k-2}{n-1}.
\]
\label{lem:GC-sum-B}
\end{lemma}
\begin{proof}
This is obtained in a similar manner to the preceding Lemma, except there is an additional box of content 0 (the one at the end of the row of length 2), and one box of content $n-k-1$ has been removed. 
\end{proof}

\begin{thm}
Let $n\geq 4$. When $2\leq p \leq \frac{n}{2}$, number of $(p,n-1,n)$-dipoles in an orientable surface of genus $g$ is
\[
(n-2)![t^{n-2g}] D_{n,p}(t),
\]
where
\begin{align*}
D_{n,p}(t) &=  \binom{t+n-1}{n} + \sum_{p\leq k\leq n-p-1} \frac{(n-1)(n-p)}{k(n-k-1)}  \binom{t+n-k-1}{n} \\
&\qquad+ \binom{t}{n} - \sum_{p-1 \leq k\leq n-p-1}  \frac{(n-p)}{(n-k-1)(k+1)} t\binom{t+n-k-2}{n-1}.
\end{align*}
When $\frac{n}{2}\leq p\leq n-1$, the generating series for $(p,n-1,n)$ dipoles is 
\begin{align*}
D_{n,p}(t) &=  \binom{t+n-1}{n} - \sum_{n-p\leq k\leq p-1} \frac{(n-1)(n-p)}{k(n-k-1)}  \binom{t+n-k-1}{n} \\
&\qquad+ \binom{t}{n} + \sum_{n-p \leq k\leq p-2}  \frac{(n-p)}{(n-k-1)(k+1)} t\binom{t+n-k-2}{n-1}.
\end{align*}
\label{thm:p1n-full-solution}
\end{thm}
\begin{proof}
By the Euler-Poincar\'{e} formula, the number of $(p,n-1,n)$-dipoles on a surface of genus $g$ may be obtained by summing the expression given in Theorem \ref{thm:p1n-dipoles-general-formula} over all partitions having $n-2g$ parts. Applying Lemmas \ref{lem:GC-sum-A} and \ref{lem:GC-sum-B} to this sum yields
\begin{align*}
D_{n,p}(t) &=  \binom{t+n-1}{n} + \sum_{1\leq k\leq n-2} (-1)^k(n-1) \left(\frac{\gamma^{(n-k,1^k),n-k}_{(p,n-p),p}}{n-k-1} - \frac{\gamma^{(n-k,1^k),1}_{(p,n-p),p}}{k} \right) \binom{t+n-k-1}{n} \\
&\qquad+ \binom{t}{n} + \sum_{1\leq k\leq n-3} (-1)^k \frac{k(n-k-2) \gamma^{(n-k-1,2,1^{k-1}),2}_{(p,n-p),p}}{(n-k-1)(k+1)} t\binom{t+n-k-2}{n-1}.
\end{align*}
This expression may be further simplified by using Lemmas \ref{cor:gen-char-two-part-evaluations-A}, \ref{lem:gen-char-two-part-evaluations-B} and \ref{lem:gen-char-two-part-evaluations-C}. We consider two cases as follows.

\noindent \textbf{Case 1}: $2\leq p \leq \frac{n}{2}$.
In this case, when $1 \leq k< p$, 
\begin{align*}
\frac{\gamma^{(n-k,1^k),n-k}_{(p,n-p),p}}{n-k-1} - \frac{\gamma^{(n-k,1^k),1}_{(p,n-p),p}}{k} &= \frac{(-1)^k}{n-1} - \frac{(-1)^k}{n-1}= 0.
\end{align*}
When $p\leq k \leq n-p-1$,
\begin{align*}
\frac{\gamma^{(n-k,1^k),n-k}_{(p,n-p),p}}{n-k-1} - \frac{\gamma^{(n-k,1^k),1}_{(p,n-p),p}}{k} &=  \frac{(-1)^k(n-p)}{(n-1)(n-k-1)} - \frac{(-1)^{k-1}(n-p)}{(n-1)k} \\
&= \frac{(-1)^k(n-p)}{k(n-k-1)}.
\end{align*}
When $n-p\leq k \leq n-2$,
\begin{align*}
\frac{\gamma^{(n-k,1^k),n-k}_{(p,n-p),p}}{n-k-1} - \frac{\gamma^{(n-k,1^k),1}_{(p,n-p),p}}{k} &= \frac{(-1)^{k-1}}{n-1} - \frac{(-1)^{k-1}}{n-1} = 0.
\end{align*}
When $k<p-1$ or $k\geq n-p$, then $\gamma^{(n-k-1,2,1^{k-1}),2}_{(p,n-p),p}=0$. For $p-1 \leq k \leq n-p-1$,
\[
\gamma^{(n-k-1,2,1^{k-1}),2}_{(n-p,p),p} = \frac{(-1)^{k-1}(n-p)}{k(n-k-2)}.
\]
Combining these facts gives the result in the statement of the Theorem. 

\noindent \textbf{Case 2}: $\frac{n}{2}\leq p\leq n-1$. In this case, when $1 \leq k< n-p$, 
\begin{align*}
\frac{\gamma^{(n-k,1^k),n-k}_{(p,n-p),p}}{n-k-1} - \frac{\gamma^{(n-k,1^k),1}_{(p,n-p),p}}{k} &= \frac{(-1)^k}{n-1} - \frac{(-1)^k}{n-1} = 0.
\end{align*}
When $n-p\leq k \leq p-1$,
\begin{align*}
\frac{\gamma^{(n-k,1^k),n-k}_{(p,n-p),p}}{n-k-1} - \frac{\gamma^{(n-k,1^k),1}_{(p,n-p),p}}{k} &=  \frac{(-1)^{k-1}(n-p)}{(n-1)(n-k-1)} - \frac{(-1)^{k}(n-p)}{(n-1)k} \\
&= \frac{(-1)^{k-1}(n-p)}{k(n-k-1)}.
\end{align*}
When $p\leq k \leq n-2$,
\begin{align*}
\frac{\gamma^{(n-k,1^k),n-k}_{(p,n-p),p}}{n-k-1} - \frac{\gamma^{(n-k,1^k),1}_{(p,n-p),p}}{k} &= \frac{(-1)^{k-1}}{n-1} - \frac{(-1)^{k-1}}{n-1} = 0.
\end{align*}
When $k<n-p$ or $k> p-2$, then $\gamma^{(n-k-1,2,1^{k-1}),2}_{(p,n-p),p}=0$. For $n-p \leq k \leq p-2$,
\[
\gamma^{(n-k-1,2,1^{k-1}),2}_{(n-p,p),p} = \frac{(-1)^{k}(n-p)}{k(n-k-2)}.
\]
Combining these facts gives the result in the statement of the Theorem. 
\end{proof}
This theorem may be used to prove the following result.
\begin{thm}
\label{thm:p1n-symmetry}
Let $p\geq 2$, and let $p,p'$ be such that $p+p'=n+1$. Then the number of $(p,n-1,n)$-dipoles on a surface of genus $g$ is equal to the number of $(p',n-1,n)$-dipoles on the same surface.
\end{thm}
\begin{proof}
We may assume that $p\leq \frac{n}{2}$, since the case $p=p'=\frac{n+1}{2}$ is trivial. By Theorem \ref{thm:p1n-full-solution}, it suffices to prove that $D_{n,p}-D_{n,p'}=0$. By routine simplification,
\begin{align*}
D_{n,p}-D_{n,p'} &=  \frac{n-1}{n-p} \binom{t+n-p}{n} + \frac{n-1}{n-p}\binom{t+p-1}{n} \\
&\qquad+ \sum_{p\leq k\leq n-p-1} \frac{(n-1)^2}{k(n-k-1)} \binom{t+n-k-1}{n} \\
&\qquad- \sum_{p-1\leq k\leq n-p-1} \frac{(n-1)}{(n-k-1)(k+1)} t \binom{t+n-k-2}{n-1}.
\end{align*}
Applying the identity
\[
t\binom{t+n-k-2}{n-1} = n \binom{t+n-k-1}{n} - (n-k-1)\binom{t+n-k-2}{n-1}
\]
and simplifying gives
\begin{align*}
\frac{D_{n,p} - D_{n,p'}}{n-1} &= -\frac{1}{p} \binom{t+n-p}{n} + \frac{1}{n-p}\binom{t+p}{n}  \\
&\qquad+ \sum_{p\leq k\leq n-p-1} \frac{1}{k(k+1)} \binom{t+n-k-1}{n} \\
&\qquad + \sum_{p\leq k\leq n-p-1} \frac{1}{k} \binom{t+n-k-1}{n-1}.
\end{align*}
Since $\binom{t+n-k-1}{n-1} = \binom{t+n-k}{n}-\binom{t+n-k-1}{n}$, then
\begin{align*}
\frac{D_{n,p} - D_{n,p'}}{n-1}&= -\frac{1}{p} \binom{t+n-p}{n} + \frac{1}{n-p}\binom{t+p}{n}  \\
&\qquad+ \sum_{p\leq k\leq n-p-1} \left(\frac{1}{k} \binom{t+n-k}{n} - \frac{1}{k+1}\binom{t+n-k-1}{n} \right)=0 .
\end{align*}
\end{proof}
The question of finding a combinatorial proof of this result is an open problem. The algebraic methods employed in this paper provide some hints as to what such a proof would look like. Because of Theorem \ref{thm:p1n-dipoles-general-formula} and the fact that generalized characters are linearly independent, the symmetry described in Theorem \ref{thm:p1n-dipoles-general-formula} is not present when information about face structure is retained. Thus, any bijection between $(p,n-1,n)$-dipoles and $(p',n-1,n)$-dipoles cannot preserve face degree sequences.

\subsection{Generating series for $Z_1$-decompositions of a full cycle}

We recall the problem of near-central decompositions of a full cycle, defined in Section \ref{sec:Z1-decompositions}, which we are now in a position to address.
\begin{thm}
\label{thm:Z1-cycle-decompositions}
Let $\lambda,\mu \vdash n$. Let $i$ be a part of $\lambda$ and let $j$ be a part of $\mu$. Let $C\in \mathcal{C}_{(n)}$. The number of factorizations $C= \sigma_1\sigma_2$ such that $\sigma_1\in \mathcal{C}_{\lambda,i}$ and $\sigma_2 \in \mathcal{C}_{\mu,j}$ is given by
\[
\frac{ |\mathcal{C}_{\lambda,i}||\mathcal{C}_{\mu,j}| }{(n-1)^2n!} \sum_{1\leq k\leq n-1} \frac{(-1)^{k-1}}{\binom{n-2}{k-1}} [x^ky^k] T_{n,i,j}(x,y) H_{\lambda \setminus i}(x) H_{\mu \setminus j}(y),
\]
where 
\[
T_{n,i,j} = xyR_{n,i}(x)R_{n,j}(y) - S_{n,i}(x)S_{n,j}(y),
\]
and $R_{n,i}$ and $S_{n,i}$ are defined in Lemma \ref{lem:hook-genchar-series}.
\end{thm}
\begin{proof}
The number of such factorizations is given by
\begin{align*}
[K_{(n),n}] K_{\lambda,i}K_{\mu,j} &= c^{(n),n}_{\lambda,i,\mu,j} \\
&=  \frac{ |\mathcal{C}_{\lambda,i}||\mathcal{C}_{\mu,j}| }{n!} \sum_{\rho \vdash n, \ell \in \rho} \frac{\gamma^{\rho,\ell}_{\lambda,i} \gamma^{\rho,\ell}_{\mu,j}\gamma^{\rho,\ell}_{(n),n}}{d_{\ell_-(\rho)}} \frac{d_{\rho}}{d_{\ell_-(\rho)}}.
\end{align*}
It can be shown, in a manner similar to the proof of Lemma \ref{lem:GC-evaluation-(n-1)-cycle}, that
\[
\gamma^{\mu,j}_{(n),n} = 
\begin{cases}
(-1)^{k}\frac{n-k-1}{n-1} & \text{~if } \mu = (n-k,1^{k}), j=n-k; \\
(-1)^{k}\frac{k}{n-1} &\text{~if } \mu = (n-k,1^{k}), j=1; \\
0 & \text{~otherwise.}
\end{cases}
\]
Thus,
\begin{align*}
[K_{(n),n}] K_{\lambda,i}K_{\mu,j} &= \frac{ |\mathcal{C}_{\lambda,i}||\mathcal{C}_{\mu,j}| }{n!} \left(\sum_{0\leq k\leq n-2} \frac{\gamma^{(n-k,1^k),n-k}_{\lambda,i} \gamma^{(n-k,1^k),n-k}_{\mu,j} (-1)^k }{\binom{n-2}{k}}      \right. \\
&\qquad\qquad\qquad \left. + \sum_{1\leq k\leq n-1} \frac{\gamma^{(n-k,1^k),1}_{\lambda,i} \gamma^{(n-k,1^k),1}_{\mu,j} (-1)^k }{\binom{n-2}{k-1}}   \right).
\end{align*}
Using Lemma \ref{lem:hook-genchar-series}, 
\begin{align*}
[K_{(n),n}] K_{\lambda,i}K_{\mu,j} &= \frac{ |\mathcal{C}_{\lambda,i}||\mathcal{C}_{\mu,j}| }{(n-1)^2n!} \left( \sum_{0\leq k \leq n-2} \frac{(-1)^k}{\binom{n-2}{k}} [x^ky^k] R_{n,i}(x)H_{\lambda\setminus i}(x)R_{n,j}(y)H_{\mu \setminus j}(y) \right. \\
&\qquad\qquad\qquad \left. + \sum_{1\leq k\leq n-1} \frac{(-1)^k}{\binom{n-2}{k-1}} [x^ky^k] S_{n,i}(x) H_{\lambda \setminus i}(x) S_{n,j}(y) H_{\mu\setminus j}(y)\right),
\end{align*}
from which the result follows.
\end{proof}

\section{Concluding Comments}

We anticipate that the study of $Z_1(n)$ will lead to insight about $Z_2(n)$, the algebra needed to give a full solution to the general $(p,q,n)$-dipole problem on all orientable surfaces. A recursive solution to the $(p,q,n)$-dipole problem may be obtained by a differential approach using a Join-Cut analysis. The interested reader is directed to \cite{JacksonSloss:2011a}.

\section{References}

\bibliographystyle{amsplain}
\bibliography{near-central-2}

\end{document}